\documentclass{llncs}

\usepackage{times,epsfig,amssymb,latexsym,lscape}

\newcommand{\allEq}{\hbox{{\it allEq }}}
\newcommand{\ex}{\exists}

\newcommand{\myVc}[3]{\ensuremath{#1_#2,\ldots,#1_{#3}}}
\newcommand{\myVec}[2]{\ensuremath{#1_1,\ldots,#1_{#2}}}

\title{Propagation by Selective Initialization\\
       and Its Application to\\
       Numerical Constraint Satisfaction Problems}
\author{M.H. van Emden \and B. Moa}

\institute{Department of Computer Science, \\University of Victoria,
Victoria, Canada \\
\email{\{vanemden, bmoa\}@cs.uvic.ca},\\ WWW home page:
\texttt{http://www.cs.uvic.ca/\homedir vanemden/ }\\
}

\begin{document}
\maketitle

\begin{abstract}
Numerical analysis has no satisfactory method for the more realistic
optimization models. However, with constraint programming one can
compute a cover for the solution set to arbitrarily close
approximation. Because the use of constraint propagation for composite
arithmetic expressions is computationally expensive, consistency is
computed with interval arithmetic.  In this paper we present theorems
that support, selective initialization, a simple modification of
constraint propagation that allows composite arithmetic expressions to
be handled efficiently.

\end{abstract}

\section{Introduction}
The following attributes all make an optimization problem more
difficult:  having an objective function with an unknown and possibly
large number of local minima, being constrained, having nonlinear
constraints, having inequality constraints, having both
discrete and continuous variables.  Unfortunately, faithfully
modeling an application tends to introduce many of these attributes.
As a result, optimization problems are usually linearized, discretized,
relaxed, or otherwise modified to make them feasible according to
conventional methods.

One of the most exciting prospects of constraint programming is
that such difficult optimization problems can be solved without
these
possibly invalidating modifications.  Moreover, constraint programming
solutions are of known quality:  they yield intervals
guaranteed to contain all solutions. Equally important, constraint
programming can prove the absence of solutions.

In this paper we only consider the core of the constraint programming
approach to optimization, which is to solve a system of nonlinear
inequalities:

\begin{equation}
\label{nonLinSys}
\begin{array}{ccccccccc}
g_1(x_1&,& x_2 &,& \ldots &,& x_m)  & \leq & 0     \\
g_2(x_1 &,& x_2 &,& \ldots &,& x_m) & \leq & 0     \\
\multicolumn{9}{c}{\dotfill}            \\
g_k(x_1 &,& x_2 &,& \ldots &,& x_m) & \leq & 0     \\
\end{array}
\end{equation}
It is understood that it may happen that $g_i= -g_j$ for some pairs
$i$ and $j$, so that equalities are a special case. If this occurs,
then certain obvious optimizations are possible in the methods
described here.

The ability to solve systems such as (\ref{nonLinSys})
supports optimization in more ways
than one.  In the first place, these systems occur as conditions
in some constrained optimized problems. Moreover, one of \myVec{g}{k}
could be defined as $f(\myVec{x}{m}) - c$, where $f$ is the objective
function and where $c$ is a constant.  By repeatedly solving such a
system for suitably chosen $c$, one can find the greatest value of $c$
for which (\ref{nonLinSys}) is found to have no solution.  That value
is a lower bound for the global minimum \cite{vnmdn03}.

This approach handles nonlinear inequalities with real variables.
It also allows some or all variables to be integer by regarding
integrality as a constraint on a real variable \cite{bnldr97}.

All constraint programming work in this direction has been based on
interval arithmetic. The earliest work \cite{BNR88} used a generic
propagation algorithm based directly on
domain reduction operators for primitive arithmetic
constraints. These constraints included
$sum(x,y,z)$ defined as $x+y=z$ for all reals
$x$,
$y$, and
$z$.
Also included was
$prod(x,y,z)$ defined as $xy=z$ for all reals
$x$,
$y$, and
$z$.

This was criticized in \cite{bmcv94} which advocated the use of
composite arithmetic expression directly rather than reducing them
to primitive arithmetic constraints.
In \cite{bggp99,vnmdn01b} it was acknowledged that
the generic propagation algorithm is not satisfactory for CSPs that
derive from composite arithmetic expressions.  These papers describe
propagation algorithms that exploit the structure of such expressions
and thereby improve on what is attainable by evaluating such
expressions in interval arithmetic.
Selective Initialization was first described in \cite{vnmdn03a}.
This was done under the tacit assumption that all default domains are
$[-\infty,+\infty]$. As a result some of the theorems in that paper are
not as widely applicable as claimed.

All of these researches are motivated by the severe difficulties
experienced by conventional numerical analysis to solve practical
optimization problems.  They can be regarded as attempts to fully
exploit the potential of interval arithmetic. In this paper we also
take this point of view. We show that, though Equation
(\ref{nonLinSys}) can contain arbitrarily large expressions, only a
small modification of the generic  propagation algorithm is needed to
optimally exploit the structure of these expressions.  This is made
possible by a new canonical form for Equation (\ref{nonLinSys}) that we
introduce in this paper.  In addition to supporting our application of
constraint processing to solving systems similar to Equation
(\ref{nonLinSys}), this canonical form exploits the potential for
parallelism in such systems.

\section{A software architecture for optimization problems}
\begin{figure}[!htb]
\begin{center}


\epsfxsize=4in
\leavevmode
\epsfbox{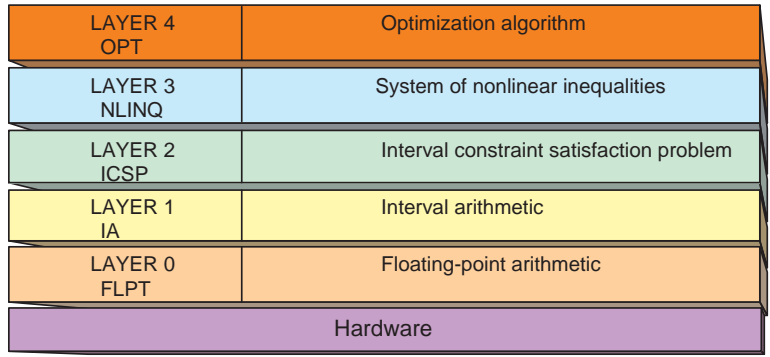}

\end{center}
\caption{
\label{softArch}
A software architecture for optimization problems.
}
\end{figure}

In Figure~\ref{softArch} we propose a hierarchical software architecture for
optimization problems. Each layer is implemented in terms of the layer below.

In the introduction we briefly remarked on how layer 4 can be reduced to layer
3. More detail is given in \cite{vnmdn03}. For the transition between layers
0 and 1 there is much material in the interval arithmetic literature. The
part that is relevant to constraint processing can be found in
\cite{hckvnmdn01}.
In the present paper we present a new method for implementing layer
3 in terms of layer 2. But first we review the transition between layers
1 and 2.

\section{Preliminaries}
In this section we provide background by reviewing some basic
concepts.  These reviews also serve to establish the terminology and
notation used in this paper.  The first few sections apply to all
constraint satisfaction problems, not only to numerical ones.

\subsection{Constraint satisfaction problems}
A {\em constraint satisfaction problem (CSP)} consists of a set of {\em
constraints}. Each of the variables in the constraint is associated
with a {\em domain}, which is the set of values that are possible for
the variable concerned. Typically, not all sets of values can be
domains. For example, sets of real values are restricted to intervals,
as described later.

A {\em valuation} is a tuple indexed by variables where the component
indexed by $v$ is an element of the domain of $v$.  A {\em solution}
is a valuation such that each constraint is true if every variable
in the constraint is substituted by the component of the valuation
indexed by the variable.  The set of solutions is a set of valuations;
hence a set of tuples; hence a \emph{relation}. We regard this relation
as the relation defined by the CSP.  In this way the relation that
is the meaning of a constraint in one CSP can be defined by
another.  This gives CSPs an hierarchical structure.

With each constraint, there is an associated {\em domain
reduction operator}\/; DRO for short.  This operator may remove
from the domains of each of the variables in the constraint certain
values that do not satisfy the constraint, given that the other
variables of the constraint are restricted to their associated
domains.  Any DRO is contracting, monotonic, and idempotent.

When the DROs of the constraints are applied in a ``fair'' order,
the domains converge to a limit or one of the domains becomes empty.
A sequence of DROs activations is \emph{fair} if every one of them occurs
an infinite number of times.  The resulting
Cartesian product of the domains becomes the greatest common fixpoint
of the DROs \cite{vnmd97,aptEssence}.  If one of the domains becomes
empty, it follows that no solutions exist within the initial domains.

In practice, we are restricted to the domains that are representable
in a computer. As there are only a finite number of these, any fair
sequence of DRO applications yields domains that remain constant
from a certain point onwards.

\subsection{Constraints}
According to the usual terminology in constraint programming,
a constraint states that a certain relation holds between
its arguments.
But in first-order predicate logic the same role
is played by an \emph{atomic formula}.
In this paper we adopt the
terminology of first-order predicate logic for the meaning
of ``atomic formula'' and we reserve ``constraint'' for a
special case.

Thus an atomic formula consists of a predicate symbol
with terms as arguments. A \emph{term} is a function
symbol with terms as arguments. What makes an atomic
formula first-order is that a predicate symbol can only
occur as the outermost symbol in the formula.

At first sight, the inequalities in Equation (\ref{nonLinSys})
are atomic formulas. In fact, they follow the related, but
different, usage that prevails in informal mathematics.
The inequality

\begin{equation}
\label{inEq}
\begin{array}{ccccccccc}
g_i(x_1&,& x_2 &,& \ldots &,& x_m)  & \leq & 0
\end{array}
\end{equation}
is an atomic formula where $\leq$ is the predicate symbol
with two arguments, which are the terms
$g_i(x_1, x_2 , \ldots , x_m)$ and $0$.
A possible source of confusion is that in mathematics
$g_i(x_1, x_2 , \ldots , x_m)$ is not necessarily interpreted
as a syntactically concrete term,
but as an abstractly conceived function $g_i$ applied to the
arguments
\myVc{x}{1}{m}. The function could be defined by means of a term
that looks quite different; such a term could be
nested and contain several function symbols. For example,
Equation (\ref{inEq}) could be
$g_i(x_1,x_2) \leq 0$
with $g_i$ defined as
$g_i(x,y) =  x^2 + xy - y^2 $ for all $x$ and $y$.
Accordingly, the atomic formula corresponding to Equation (\ref{inEq})
is
\begin{equation}
\label{inEqEx}
\begin{array}{ccccccccc}
\leq(+(sq(x),-(\times(x,y),sq(y))),0).
\end{array}
\end{equation}
Taking advantage of infix and postfix notation this becomes
$x^2 + xy - y^2 \leq 0$.

A \emph{constraint} is an atomic formula without function symbols.
An example of such an atomic formula is
$sum(x,y,z)$, which is a ternary constraint whose relation is defined
by $x+y=z$ for all reals $x$, $y$, and $z$.

In this paper we translate Equation~(\ref{inEqEx}) to a CSP
with the set of constraints
$$
\{t_1 \leq 0,sum(t_2,t_3,t_1),sq(x,t_2),sum(t_5,t_3,t_4),
  prod(x,y,t_5),sq(y,t_4)
\}
$$

Consider a constraint $c(\myVec{x}{n})$. The meaning of predicate
symbol $c$ is a relation, say $r$.
For all $i \in \{1,\ldots,n\}$, a value $a_i$ for variable
$x_i$ is \emph{inconsistent} with respect to $r$ and domains
$X_1,\ldots,X_{i-1},X_{i+1},\ldots,X_{n}$
iff it is not the case that
\begin{eqnarray*}
  \ex a_1 \in X_1,\ldots,
  \ex a_{i-1} \in X_{i-1},\ex a_{i+1}\in X_{i+1},\ldots,
  \ex a_{n} \in X_{n} \hbox{ such that } \\
 \langle \myVec{a}{n} \rangle \in r.
\end{eqnarray*}
A DRO for $c$ may replace one or more
of the domains \myVec{X}{n} by a subset of it if the set difference between
the old and the new domains contains inconsistent values only.
A DRO is \emph{optimal}
if every domain is replaced by the smallest domain containing
all its consistent values.
We call a constraint
\emph{primitive} if an optimal DRO is available for it that is
sufficiently efficiently computed.
What is sufficient depends on the context.

\subsection{Constraint propagation}
To gain information about the solution set,
inconsistent values are removed as much as possible with modest
computation effort.  For example, DROs can be applied as long as
they remove inconsistent values.  It is the task of a constraint
propagation algorithm to reach as quickly as possible a set of
domains that cannot be reduced by any DRO.
Many versions of this algorithm exist
\cite{aptEssence,rthpra01}.
They can be regarded as refinements of
the algorithm in Figure~\ref{lGPA}, which we refer to as the
\emph{generic propagation algorithm}\/; GPA for short.

GPA maintains a pool of DROs, called {\em active set}.  No order
is specified for applying these operators.  Notice that the active
set $A$ is initialized to contain all constraints.

\begin{figure}
\begin{tabbing}
put all constraints into the active set $A$\\
while \=( $A \neq \emptyset$) $\{$\\
         \>choose a constraint $C$ from $A$\\
         \>apply the DRO associated with $C$\\
         \> if one of the domains has become empty, then stop\\
         \>add to $A$ all constraints involving variables whose
           domains have changed, if any\\
         \>remove $C$ from $A$\\
$\}$
\end{tabbing}
\caption{
\label{lGPA}
A pseudo-code for GPA.
}
\end{figure}

\subsection{Intervals}
A {\em floating-point number} is any element of $F \cup \{-\infty,
+\infty\}$, where $F$ is a finite set of reals. A {\em floating-point
interval} is a closed connected set of reals, where the bounds, in
so far as they exist, are floating-point numbers. When we write
``interval'' without qualification, we mean a floating-point
interval.  For every bounded non-empty interval $X$, $lb(X)$ and
$rb(X)$ denote the least and the greatest element of $X$ respectively.
They are referred to as the left and the right bound of $X$.  If
$X$ is not bounded from below, then $lb(X) = -\infty$.
Similarly, if $X$ is not bounded from above, then
$rb(X) = +\infty$.
Thus, $ X = [lb(X), rb(X)]$ is a convenient
notation for all non-empty intervals, bounded or not.

A {\em box} is a Cartesian product of intervals.

\subsection{Solving inequalities in interval arithmetic}

Moore's idea of solving inequalities such as those in Equation
(\ref{nonLinSys}) by means of interval arithmetic is at least as important as
the subsequent applications of interval constraints to this problem.

Suppose we wish to investigate the presence of solutions of a single
inequality in Equation (\ref{nonLinSys}) in a box
$X_1 \times \cdots \times X_m$.
Then one evaluates in interval arithmetic the
expression in the left-hand side.
As values for the variables $x_1 , \ldots , x_m$
one uses the intervals $X_1 , \cdots , X_m$.
Suppose the result is the interval
$[a_i,b_i]$. We have exactly one of the following three cases.
If $0 < a_i$ for at least one $i$, then there are no solutions.
If $b_i \leq 0$ for all $i$, then all tuples
in $X_1 \times \cdots \times X_m$ are solutions.
If $a_i \leq 0 < b_i$ for at least one $i$, then either of the above may be
true.
In this same case of $a_i \leq 0 < b_i$,
it may also be that some of the tuples
in $X_1 \times \cdots \times X_m$ are solutions,
while others are not.

Again, in the case where $a_i \leq 0 < b_i$,
it may be possible to split $X_1 \times \cdots \times X_m$.
In this way, a more informative outcome may be obtained
for one or both of the results of splitting.
Such splits can be repeated as long as possible and necessary.


\subsection{Interval CSPs}

Problems in a wide variety of application areas can be expressed
as CSPs.  Domains can be as different as booleans, integers, finite
symbolic domains, and reals.  In this paper we consider {\em Interval
CSPs} (ICSPs), which are CSPs where all domains are intervals and
all constraints are primitive.

ICSPs are important because they encapsulate what can be
efficiently computed; they represent Layer 2 in the software
architecture of Figure~\ref{softArch}.
The layer is distinct from Layer 3 because in Equation (\ref{nonLinSys})
there typically occur atomic formulas
that contain function symbols.

To emphasize the role of ICSPs as a layer of software
architecture, we view them as a \emph{virtual machine}, with a
function that is similar to those for Prolog or Java. Just as a
program in Prolog or Java is translated to virtual machine
instructions, a system such as Equation (\ref{nonLinSys}) can be
translated to an ICSP, as described in a later section.

The instructions of the ICSP level are DROs, one for each
constraint.  As an example of such an ICSP virtual machine instruction,
consider the DRO for product constraint. It reduces the box
$[a,b]\times [c,d]\times [e,f]$ to the box that has the projections
\begin{eqnarray}
\varphi([a,b] & \cap & ([e,f]/[c,d])) \nonumber\\
\varphi([c,d] & \cap & ([e,f]/[a,b])) \nonumber\\
\varphi([e,f] & \cap & ([a,b]*[c,d])) \nonumber
\end{eqnarray}
Here $\varphi$ is the function that yields the smallest interval
containing its argument.

Of particular interest is the effect of the DRO when all variables
have $[-\infty,+\infty]$ as domain.  For each argument, the domain
after application of the DRO is defined as the \emph{default domain}
of that argument.  Typically, default domains are $[-\infty,+\infty]$.
Notable exceptions include the constraint $sin(x,y)$ (defined as
$y = \sin(x)$), where the default domain of $y$ is $[-1,1]$. Another
is $\hbox{sq}(x,y)$ (defined as $y = x^2$), where the default domain
of $y$ is $[0,\infty]$.

A difference with other virtual machines is that
a program for the ICSP virtual machine is an unordered collection of DROs.
Programs for other virtual machines are ordered sequences of instructions.
In those other virtual machines,
the typical instruction does not specify the successor instruction.
By default this is taken to be the next one in textual order.
Execution of the successor is implemented
by incrementing the instruction counter by one.

The simplicity of the instruction sequencing in conventional virtual
(and actual) machines is misleading. Many instruction executions
concern untypical instructions, where the next instruction is
specified to be another than the default next instruction.  Examples
of such untypical instructions are branches (conditional or
unconditional) and subroutine jumps.

In the ICSP virtual machine, the DROs are the instructions,
and they form an unordered set. Instead of an instruction counter
specifying the next instruction, there is the active set of GPA
containing the set of possible next instructions. Instead of an
instruction or a default rule determining the next instruction to
be executed, GPA selects in an unspecified way which of the DROs
in the active set to execute.  In this way, programs can be
declarative:  instructions have only meaning in terms of \emph{what}
is to be computed.  \emph{How} it is done (instruction sequencing),
is the exclusive task of the virtual machine.

\subsection{A canonical form for nonlinear numerical inequalities}
\label{arch}

Equation (\ref{nonLinSys}) may have multiple occurrences of variables
in the same formula.  As there are certain advantages in avoiding such
occurrences, we rewrite without loss of generality the system in
Equation~(\ref{nonLinSys}) to the canonical form shown in
Figure~\ref{singleSys}.

\begin{figure}
\begin{displaymath}
\begin{array}{ccccccccc}
g_1(y_1&,& y_2 &,& \ldots &,& y_n) & \leq &0     \\
g_2(y_1 &,& y_2 &,& \ldots &,& y_n) & \leq & 0     \\
\multicolumn{9}{c}{\dotfill}            \\
g_k(y_1 &,& y_2 &,& \ldots &,& y_n) & \leq & 0     \\

\allEq(v_{1,1}&,& v_{1,2}&,& \ldots&,& v_{1,{n_1}}) &&  \\
\multicolumn{7}{c}{\dotfill}            \\
\allEq(v_{p,1}&,& v_{p,2}&,& \ldots&,& v_{p,{n_p}}) &&  \\
\end{array}
\end{displaymath}
\caption{
\label{singleSys}
A system of non-linear inequalities without multiple occurrences
of variables. Instead,
the set $\{y_1,\ldots,y_n\}$ is partitioned into equivalence
classes
$V_1,\ldots,V_p$ where
$V_j$ is a subset $\{v_{j,1}, \ldots, v_{j,{n_j}}\}$
of $\{y_1,\ldots,y_n\}$,
for $j \in \{1,\ldots, p\}$.
An $\mathit{allEq}$ constraint asserts that its arguments
are equal.
}
\end{figure}

In Figure~\ref{singleSys}, the expressions for the functions $g_1,
\ldots, g_k$ have no multiple occurrences of variables.  As a
result, they have variables \myVec{y}{n} instead of \myVec{x}{m},
with $m \leq n$ as in Equation~(\ref{nonLinSys}).  This canonical
form is obtained by associating with each of the variables $x_j$
in Equation~(\ref{nonLinSys}) an equivalence class of the variables
in Figure~\ref{singleSys}.  This is done by replacing in
Equation~(\ref{nonLinSys}) each occurrence of a variable by a
different element of the corresponding equivalence class.
This is possible by ensuring
that each equivalence class is as large as the largest number of
multiple occurrences.  The predicate $\allEq$ is true if and only
if all its real-valued arguments are equal.

An advantage of this translation is that evaluation in interval
arithmetic of each expression gives the best possible result, namely
the range of the function values. At the same time, the $\allEq$
constraint is easy to enforce by making all intervals of the
variables in the constraint equal to their common intersection.
This takes information into account from all $k$ inequalities.  If
the system in its original form as in Equation~\ref{nonLinSys},
with multiple occurrences, would be translated to a CSP, then only
multiple occurrences in a single expression would be exploited at
one time.

In the coming sections and without loss of generality, we will only
consider expressions without multiple occurrences of variables.

\subsection{Translating nonlinear inequalities to ICSPs}
\label{translation}
ICSPs represent what we \emph{can} solve.  They consist of atomic
formulas without function symbols that, moreover, have efficient
DROs.  Equation (\ref{nonLinSys}) exemplifies what we \emph{want}
to solve: it consists of atomic formulas typically containing deeply
nested terms.

\paragraph{The tree form of a formula}
We regard a first-order atomic formula as a tree.  The unique
predicate symbol is the root.  The terms that are the arguments of
the formula are also trees and they are the subtrees of the root.
If the term is a variable, then the tree only has a root, which is
that variable.  A term may also be a function symbol with one or
more arguments, which are terms. In that case, the function symbol
is the root with the argument terms as subtrees.

In the tree form of a formula the leaves are variables.
In addition,
we label every node that is a function symbol with a unique variable.

Any constants that may occur in the formula are replaced by unique
variables.  We ensure that the associated domains contain the
constants and are as small as possible.

\paragraph{Translating a formula to an ICSP}
The tree form of a formula thus labeled is readily translated to an
ICSP. The translation has a set of constraints in which each
element is obtained by translating an internal node of the tree.

The root translates to $p(\myVc{x}{0}{n-1})$
where $p$ is the predicate symbol that is the root and
\myVc{x}{0}{n-1}
are the variables labeling the children of the root.

A non-root internal node
of the form $f(\myVc{t}{0}{n-1})$
translates to
$F(\myVc{x}{0}{n-1}, y)$, where
\begin{itemize}
\item
$y$ is the variable labeling the node
\item
\myVc{x}{0}{n-1} are the variables labeling the child nodes
\item
$F$ is the relation defined by $F(\myVc{a}{0}{n-1}, v)$
iff $v= f(\myVc{a}{0}{n-1})$ for all $\myVc{a}{0}{n-1},v$.
\end{itemize}

\subsection{Search}

Propagation may terminate with small intervals for all variables
of interest. This is rare. More likely, propagation leaves a large
box as containing all solutions, if any. To obtain more information
about possibly existing solutions, it is necessary to split an ICSP into two
ICSPs and to apply propagation to both. An ICSP $S'$ is a result
of splitting an ICSP $S$ if $S'$ has the same constraints as $S$
and differs only in the domain for one variable, say, $x$.  The
domain for $x$ in $S'$ is the left or right half of the domain for
$x$ in $S$.

A \emph{search strategy} for an ICSP is a binary tree representing
the result of successive splits. Search strategies can differ
greatly in the effort required to carry them to completion.  The
most obvious search strategy is the \emph{greedy strategy}\/: the
one that ensures that all intervals become small enough by choosing
a widest domain as the one to be split.

This is a plausible strategy in the case where the ICSP has a few
point solutions. In general, the set of solutions is a continuum:
a line segment, a piece of a surface, or a variety in a higher
dimensional space that has positive volume.  In such cases we prefer
the search to result in a single box containing all solutions.
Of course we also prefer such a box to be as small as possible. The
greedy search strategy splits the continuum of solutions into an
unmanageably large number of small boxes.  It is not clear that
the greedy strategy is preferable even in the case of a few
well-isolated point solutions. In general we need a search
strategy other than the greedy one.

A more promising search strategy was first described in the \verb+absolve+
predicate of the BNR Prolog system \cite{BNR88} and by \cite{bmcv94},
where it is called \emph{box consistency}.

The box consistency search strategy selects a variable and a domain
bound. Box consistency uses binary search to
determine a boundary interval that can be shown to contain no
solutions. This boundary interval can then be removed from the
domain, thus shrinking the domain.  This is repeated until a boundary
interval with width less than a certain tolerance is found that
cannot be shown to contain no solutions. When this is the case for
both boundaries of all variables, the domains are said to be \emph{box
consistent} with respect to the tolerance used and with respect to
the method for showing inconsistency.
When this method is interval arithmetic, we obtain
\emph{functional box consistency}. When it is
propagation, then it is called \emph{relational box consistency}
\cite{vnmdn01b}.

All we need to know about search in this paper is that greedy search
and box consistency are both search strategies and that both can
be based on propagation.  Box consistency is the more promising
search strategy. Thus we need to compare interval arithmetic and
propagation as ways of showing that a nonlinear inequality has no
solutions in a given box.
This we do in section~\ref{PSI}.

\section{Propagation with selective initialization}
\label{PSI}

Suppose we have a term that can be evaluated in interval
arithmetic. Let us compare the interval that is the result of such
an evaluation with the effect of GPA on the ICSP obtained by
translating the term as described in section~\ref{translation}.

To make the comparison possible we define evaluation of a term
in interval arithmetic.
The definition follows the recursive structure of the term:
a term is either a variable or
it is a function symbol with terms as arguments.
If the term is an interval, then the result is that interval.
If the argument is function $f$ applied to arguments, then the result
is $f$ evaluated in interval arithmetic applied to the results of
evaluating the arguments in interval arithmetic.
This assumes that every function symbol denotes a function that is
defined on reals as well as on intervals.
The latter is called the \emph{interval extension} of the former.
For a full treatment of interval extensions,
see \cite{moore66,nmr90,hnsn92}.
The following lemma appears substantially as Theorem 2.6 in \cite{chn98}.

\begin{lemma}
\label{basic}
Let $t$ be a term
that can be evaluated in interval arithmetic.
Let the variables of $t$ be $x_1,\ldots,x_n$.
Let $y$ be the variable associated with the root of the tree form
of $t$.
Let $S$ be the ICSP that results from translating $t$, where the domains
of $x_1,\ldots,x_n$ are $X_1,\ldots,X_n$ and where the domains
of the internal variables are $[-\infty,+\infty]$.
After termination of GPA applied to $S$, the domain for $y$
is the interval that results from interval arithmetic evaluation of $t$.
\end{lemma}

\begin{proof}
Suppose that a variable of a constraint has domain $[-\infty,+\infty]$.
After applying the DRO for that constraint, this domain has become
the result of the interval arithmetic operation that obtains the
domain for this variable from the domains of the other variables
of the constraint.\\
According to \cite{vnmd97,aptEssence}, every fair sequence of DROs
converges to the same domains for the variables.  These are also
the domains on termination of GPA.  Let us consider a fair sequence
that begins with a sequence $s$ of DROs that mimics the evaluation
of $t$ in interval arithmetic. At the end of this, $y$ has the
value computed by interval arithmetic.  This shows that GPA gives
a result that is a subinterval of the result obtained
by interval arithmetic.\\
GPA terminates after activating the DROs in $s$. This is because
in the interval arithmetic evaluation of $t$ an operation is only
performed when its arguments have been evaluated. This means that
the corresponding DRO only changes one domain. This domain is the
domain of a unique variable that occurs in only one constraint that
is already in the active set. Therefore none of the DRO activations
adds a constraint to the active set, which is empty after $s$.
\end{proof}

GPA yields the same result whatever the way constraints are selected
in the active set.  Therefore GPA always gives the result of interval
arithmetic evaluation. However, GPA may obtain this result in an
inefficient way by selecting constraints that have no effect.
This suggests that the active set be structured in a way that reflects
the structure of $t$. This approach has been taken in \cite{bggp99,vnmdn01b}.

The proof shows that, if the active set had not contained any of
the constraints only involving internal variables, these constraints
would have been added to the active set by GPA. This is the main
idea of selective initialization.
By initializing and ordering the active set
in a suitable way and leaving GPA otherwise unchanged,
it will obtain the interval arithmetic result
with no more operations than interval arithmetic.
This assumes the optimization implied by the Totality Theorem in \cite{hckmdw98}.


\begin{definition}
A constraint is a \emph{seed constraint}
iff at least one of its variables has a domain that differs from the default
domain assigned to that variable.
\end{definition}

For example, the term $\sin(x_1)+\sin(x_2)$ translates to an ICSP
with constraints $sum(u,v,y)$, $sin(x_1,u)$, and $sin(x_2,v)$.
When the domains are $[-\infty,+\infty]$ for all variables, then the
seed constraints are $sin(x_1,u)$ and $sin(x_2,v)$.
When the domains are $[-\infty,+\infty]$ for
$x_1$, $x_2$, and $y$;
$[-1,1]$ for $u$ and $v$, then $sum(u,v,y)$ is the
one seed constraint.
When the domains are $[-\infty,+\infty]$ for
$x_2$ and $y$;
$[-1,1]$ for $x_1$, $u$ and $v$, then the seed constraints are $sum(u,v,y)$ and
$sin(x_1,u)$.

\begin{definition}
Let PSI (Propagation with Selective Initialization) be GPA except for
the following modifications.\\
(a) PSI only applies to ICSPs generated by translation from
an atomic formula.  \\
(b) The active set is a priority queue that is ordered according to the
distance from the root of the node that generated the constraint.
The greater that distance, the earlier the item is removed from
the queue.  \\
(c) The active set contains all seed constraints and no other ones.
\end{definition}

Lemma~\ref{basic} says that GPA simulates interval arithmetic as far as the
result is concerned. It does not say anything about the efficiency with which
the result is obtained. Theorem~\ref{simulate}
says that PSI obtains the result as efficiently as it is done in interval
arithmetic. This assumes the functionality optimization in the DROs
\cite{hckmdw98}.

\begin{theorem}
\label{simulate}
Let $S$ be the ICSP obtained by translating
a term $t$ in variables
\myVc{x}{1}{n},
where these variables have domains
\myVc{X}{1}{n}.
Applying PSI to S terminates after selecting no constraint more than once.
Moreover, the root variable ends up with $Z$ as domain where $Z$ is the
interval resulting from evaluating $t$ with
\myVc{x}{1}{n} substituted by \myVc{X}{1}{n}.
\end{theorem}

\begin{proof}
Suppose GPA is applied to $S$ in such a way that all non-seed
constraints are selected first.  The execution of the DRO corresponding
to these non-seed constraints does not change any domains.  Therefore
these DRO executions do not add any constraints.  As a result, the
effect of applying GPA is the same as when the active set would
have been initialized with only the seed constraints.\\
Suppose the seed constraints are selected according to priority
order.  This ensures that no future constraint selection re-introduces
a constraint previously selected. Thus GPA terminates after activating
every seed constraint exactly once.\\
Such an execution of GPA coincides step by step with that of PSI.
As GPA terminates with the correct result, so does PSI.
\end{proof}

\section{Using ICSPs to solve inequalities}

We briefly reviewed how interval arithmetic can solve systems of nonlinear
inequalities. The fundamental capability turned out to be that of evaluating
a term in interval arithmetic.
We saw that this can also be done by applying propagation to ICSPs
generated from arithmetic terms.
We now investigate how to extend this
to the use of ICSPs to solve nonlinear inequalities.

\subsection{Using ICSPs to solve a single inequality}

Suppose $S$ is the ICSP resulting from translating
$$
g_i(\myVc{x}{1}{m}) \leq 0
$$
Let $y$ be the variable labeling the left child of the root; that is,
the variable representing the value of the left-hand side.
Let \myVc{X}{1}{m} be the domains in $S$ of  \myVc{x}{1}{m}, respectively.

Now suppose that GPA is applied to $S$. One possible initial sequence
of DRO activations is the equivalent of interval arithmetic evaluation
of the left-hand side, leaving $y \leq 0$ as the only constraint in the
active set with the domain for $y$ equal to $[a_i,b_i]$, the value in
interval arithmetic of $g_i(\myVc{X}{1}{m})$.

At this stage the DRO for $y \leq 0$ is executed.
If $0 < a_i$, then failure occurs.
If $b_i \leq 0$, then the domain for $y$ is unchanged. Therefore,
no constraint is added to the active set. Termination occurs with
nonfailure. There is no change in the domain of any of
\myVc{x}{1}{m}.
The third possibility is that $a_i \leq 0 < b_i$.
In this case, the domain for $y$ shrinks: the upper bound decreases from
$b_i$ to $0$.
This causes the constraints to be brought into the active set
that correspond to nodes at the next lower level in the tree.
This propagation may continue all the way down to the lowest level
in the tree, resulting in shrinking of the domain of one or more of
\myVc{x}{1}{m}.

Let us compare this behaviour with the use of interval arithmetic to solve
the same inequality.
In all three cases, GPA gives the same outcome as interval arithmetic:
failure or nonfailure. In the first two cases, GPA gives no more information
than interval arithmetic. It also does no more work.

In the third case,
GPA may give more information than interval arithmetic:
in addition to the nonfailure outcome,
it may shrink the domain of one or more of \myVc{x}{1}{m}.
This is beyond the capabilities of interval arithmetic,
which is restricted to transmit information
about arguments of a function to information
about the value of the function.
It cannot transmit information in the reverse direction.
To achieve this extra capability,
GPA needs to do more work than
the equivalent of interval arithmetic evaluation.

In the above, we have assumed that
GPA magically avoids selecting constraints in a way that is not optimal.
In such an execution of GPA we can recognize two phases:
an initial phase that
corresponds to evaluating the left-hand side in interval arithmetic,
followed by a second phase that starts with
the active set containing only the constraint $y \leq 0$.
When we consider the nodes in the tree that correspond to the constraints
that are selected, then it is natural to call the first phase bottom-up
(it starts at the leaves and ends at the root) and
the second phase top-down
(it starts at the root and may go down as far as to touch some of the leaves).
The bottom-up phase can be performed automatically by the PSI algorithm.

The start of the top-down phase is similar to
the situation that occurs in search.
In both search and in the top-down phase a different form of selective
initialization can be used, shown in the next section.

The bottom-up phase and the top-down phase are separated by a state
in which the active set only contains $y \leq 0$.
For reasons that become apparent in the next section,
we prefer a separate treatment of this constraint:
not to add it to the active set and to execute the shrinking of
the domain for $y$ as an extraneous event.
This is then a special case of termination of GPA, or its equivalent PSI,
followed by the extraneous event of shrinking one domain.

The Pseudo-code for PSI algorithm is given in Figure~\ref{PSIAlg}.

\begin{figure}[!htb]
\begin{tabbing}
let the active set $A$ be a priority queue in which the
constraints are\\
$\;\;$ ordered according to the level they occupy in the tree,\\
$\;\;$ with those that are further away from the root
placed nearer to the front of the queue\\
put only \textbf{seed} constraints into $A$\\
while \=( $A \neq \emptyset$) $\{$\\
         \>choose a constraint $C$ from $A$\\
         \>apply the DRO associated with $C$\\
         \> if one of the domains is empty, then stop\\
         \>add to $A$ all constraints involving variables whose
           domains have changed, if any\\
         \>remove $C$ from $A$\\
$\}$
\end{tabbing}
\caption{
\label{PSIAlg}
Pseudo-code for Propagation with Selective Initialization (PSI).
}
\end{figure}
The correctness of PSI algorithm can be easily deduced from
the following theorem.

\begin{theorem}
\label{modEval}
Consider the ICSP $\mathcal{S}$ obtained from the tree $T$ of the atomic
formula
$g_i(x_1,\ldots,x_n) \leq 0$.
Suppose we modify GPA so that the active set is initialized to contain
instead of all constraints only seed constraints.
Suppose also that the active set is a priority queue in which the
constraints are ordered according to the level they occupy in
the tree $T$, with those that are further away from the root
placed nearer to the front of the queue.  Then GPA terminates
with the same result as when the active set would have been
initialized to contain all constraints.
\end{theorem}

\begin{proof}
As we did before, suppose that in GPA the active set $A$ is initialized
with all constraints such that the seed constraints are at the end of the
active set. Applying any DRO of a constraint that is not a seed constraint
will not affect any domain. Thus, the constraints that are not seed
constraints can be removed from the active set without changing the result
of GPA. Since the GPA
does not specify any order, $A$ can be ordered as desired.
Here we choose to order it in such a way we get an efficient GPA when
used to evaluate an expression (see previous section).
\end{proof}

\section{Selective Initialization for search}

Often we find that after applying GPA to an ICSP $S$,
the domain $X$ for one of the variables, say $x$, is too wide.
Search is then necessary.
This can take the form of splitting $S$ on the domain for $x$.
The results of such a split are
two separate ICSPs $S_1$ and $S_2$
that are the same as $S$ except for the domain of $x$.
In $S_1$, $x$ has as domain the left half of $X$;
in $S_2$, it is the right half of $X$.

However, applying GPA to $S_1$ and $S_2$ entails duplicating
work already done when GPA was applied to $S$.
When splitting on $x$ after termination of the application
of GPA to $S$, we have the same situation as at the beginning
of the downward phase of applying GPA to an inequality:
the active set is empty and an extraneous event changes
the domain of one variable to a proper subset.

The following theorem justifies a form of the PSI algorithm where the
active set is initialized with what is typically a small subset of
all constraints.

\begin{theorem}
\label{modSolveGeneral}
Let $T$ be the tree obtained from the atomic formula
$g_i(x_1,\ldots,x_n) \leq 0$.
Let $\mathcal{S}$ be the ICSP obtained from $T$.
Let $x$ be a variable in $\mathcal{S}$.
Suppose we apply GPA to $\mathcal{S}$.
After the termination of GPA, suppose the domain
of $x$ is changed to an interval that is a proper subset of it.
If we apply GPA to $\mathcal{S}$ with an active set initialized
with the constraints only involving $x$,
then GPA terminates with the same result as
when the active set would have been initialized to contain
all constraints.
\end{theorem}

\begin{proof}
To prove Theorem~\ref{modSolveGeneral}, we should
show that initializing GPA with all constraints
gives the same results as when it is initialized with only the
constraints involving $x$.

Since no ordering is specified for the active set of GPA,
we choose an order in which the constraints involving $x$ are
at the end of the active set.
Because DROs are idempotent,
all constraints at the front of the active set,
different from those involving $x$, do not affect
any domain.
Thus removing them from the active set in the initialization
process does not change the fixpoint of the GPA.
Thus, Theorem~\ref{modSolveGeneral} is proved.
\end{proof}

\section{Further work}
We have only considered the application of selective initialization
to solve a single inequality.
A conjunction of inequalities such as Equation~(\ref{nonLinSys})
can be solved by solving each in turn.
This has to be iterated because the solving of another inequality
affects the domain of an inequality already solved.
This suggests performing the solving of all inequalities in parallel.
Doing so avoids the waste of completing an iteration
on the basis of unnecessarily wide intervals.
It also promises speed-up because many of the DRO activations
only involve variables that are unique to the inequality.
In the current version of the design of our algorithm,
we combine this parallelization with a method of minimizing
the complexity usually caused by multiple occurrences of variables.

\section{Conclusions}
Before interval methods it was not clear how to tackle
numerically realistic optimization models.
Only with the advent of interval arithmetic in the 1960s \cite{moore66}
one could for the first time at least say:
``If only we had so much memory and so much time,
then we could solve this problem.''

Interval arithmetic has been slow in developing.
Since the 1980s constraint programming has added fresh impetus to
interval methods.
Conjunctions of nonlinear inequalities, the basis for optimization,
can be solved both with interval
arithmetic and with constraint programming. In this paper we relate
these two approaches.

It was known that constraint propagation subsumes interval arithmetic.
It was also clear that using propagation for the special case of
interval arithmetic evaluation is wasteful. In this paper we present
an algorithm for propagation by Selective Initialization
that
ensures that propagation is as efficient
in the special case of interval arithmetic evaluation.
We also apply Selective Initialization
for search and for solving inequalities.
Preliminary results on a parallel version of the methods
presented here suggest that realistic optimization models
will soon be within reach of modest computing resources.

\section*{Acknowledgments}
We acknowledge generous support by the
University of Victoria, the Natural Science and Engineering Research
Council NSERC, the Centrum voor Wiskunde en Informatica CWI, and the
Nederlandse Organisatie voor Wetenschappelijk Onderzoek NWO.

\end{document}